\documentclass[12pt]{article}
\usepackage{latexsym, amsfonts, amsmath, amssymb}

\newcommand{\qed}{\hskip 5mm \rule{2.5mm}{2.5mm}\vskip 10pt}
\newcommand{\R}{{\mathbb R}}

\newcommand{\N}{{\mathbb N}}

\newcommand{\E}{{\mathbb E}}
\newcommand{\T}{{\mathbb T}}

\renewcommand{\S}{{\mathbb S}}

\newcommand{\proof}{\noindent{\em Proof:\ }}

\begin{document}
\newtheorem{thm}{Theorem}[section]
\newtheorem{defs}[thm]{Definition}
\newtheorem{lem}[thm]{Lemma}
\newtheorem{note}[thm]{Note}
\newtheorem{cor}[thm]{Corollary}
\newtheorem{prop}[thm]{Proposition}
\renewcommand{\theequation}{\arabic{section}.\arabic{equation}}
\newcommand{\newsection}[1]{\setcounter{equation}{0} \section{#1}}
\title{Markov Processes on Riesz Spaces
       }
\author{
 Jessica Vardy${}^\sharp$\\
 Bruce A. Watson${}^\sharp$ \footnote{ Supported by NRF grants FA2007041200006 and 
IFR2010042000057, and  in part by the Centre for Applicable Analysis and
 Number Theory.} \\ \\
 ${}^\sharp$ School of Mathematics\\
 University of the Witwatersrand\\
 Private Bag 3, P O WITS 2050, South Africa }
\maketitle
\abstract{ 
 Measure-free discrete time stochastic processes in Riesz spaces 
 were formulated and studied by Kuo, Labuschagne and Watson. 
 Aspects relating martingales, stopping times, convergence of these processes
 as well as various decomposition were considered. Here we formulate and study
 Markov processes in a measure-free Riesz space setting. 
}
\parindent=0in
\parskip=.10in
\newsection{Introduction}

Markov processes have been studied extensively since their introduction in 1906, \cite{markov}, by 
Andrey Markov. Roughly speaking a Markov process is a stochastic process with the property that, 
given the present state, the past and future states are independent.
The applications of Markov processes pervade almost all areas of science, economics and 
engineering.

The early theory allowed only discrete state spaces. It was with the measure-theoretic 
formulation of probability theory in the 1930's by A. Kolmogorov that
the general theory, \cite{kolmog}, could be developed.
In addition, Markov processes have been considered as the origin of the theory of stochastic 
processes and, as such, certainly 
deserve the attention required to give a measure-free formulation of the theory. 
Markov processes demonstrate the strong link between measure theory and probability theory which 
makes the generalization to the measure-free setting
that much more challenging and interesting.

A stochastic process is traditionally defined in terms of measurable functions where the underlying 
measure space is a probability space i.e.  the measure of the whole space is one.  
As was noted by Kuo, \cite{thesis}, the underlying order structure on the spaces of measurable 
functions plays a central role in the study of stochastic processes.

The setting in which we will pose Markov processes is that of a Riesz space 
(i.e. a vector space with an order structure that is compatible with the algebraic structure on it)
 with a weak order unit.  
The study of Markov processes in Riesz spaces gives one insight into the underlying mechanisms
 of the theory and, in addition,
unifies the development of the subject for a variety of settings: spaces of measurable functions, 
Banach lattices and ${\cal L}^p$-spaces for example,
see \cite{D-U, E-S, egghe, schaefer, stoica, stoica-2}.

Rao showed that conditional expectation operators on $L^p$ spaces can be characterized as 
positive contractive projections which leave the ${\bf 1}$ function invariant (where ${\bf 1}$ is the constant function with value
1), see \cite{rao}.  This 
formed the basis for the definition of
a conditional expectation operator on a Riesz space with weak order unit.  The conditional expectation operator 
is defined as 
a positive order continuous projection that maps weak order units to weak order units and has a Dedekind complete range 
space. A more detailed explanation is given in \cite{klw-indag}.
From this foundation, in \cite{klw-indag,klw-exp,klw-conv,klw-erg}, Kuo, Labuschagne and Watson 
developed many of the fundamental results of martingale theory in the measure free context of a 
Riesz space with weak order unit. 

Markov processes for which the state spaces may be non-separable are usually defined via conditional 
expectation operators and implicitly rely on the Radon-Nikod\'ym theorem. 
In \cite{watson} a Riesz space analogue of the And\^o-Douglas-Radon-Nikod\'ym theorem
 was given.  
Building on this framework we give here a generalization of Markov processes to a Riesz 
spaces setting.

We would like to thank the referees for their valuable recommendations.

\newsection{Preliminaries}

The reader is assumed familiar with the notation and terminology of Riesz spaces, 
for details see \cite{A-A, L-Z, zaanen} or, for the more specific aspects used here, 
see \cite{klw-indag, klw-exp}.

\begin{defs}\label{conditional-exp-riesz}
 Let $E$ be a Riesz space with a weak order unit.
 A positive order continuous projection $T:E\to E,$ with
 range ${\cal R}(T)$ a Dedekind complete Riesz subspace of $E,$ 
 is called a conditional expectation 
 if $T(e)$ is a weak order unit of $E$ for each weak order unit $e$ of $E$.
\end{defs}

In the above definition the condition that $Te$ is a weak order unit
for each weak order $e$ can be replaced with there exists a weak order unit $e$ with $e=Te$, 
to yield an equivalent definition, see \cite{klw-indag} for details.

If $E=L^p(\Omega,{\cal F},\mu), 1\le p\le \infty,$ where $\mu$ is a positive measure.  
Let $\Sigma$ be a sub-$\sigma$-algebra of $\cal{F}$ and $T= \E[\cdotp | \Sigma]$.  Then 
each $f \in E$ with $f >0$ almost everywhere and $f \  \Sigma$-measurable is a weak order
unit for $E$ with $f = Tf$.  In this case, $\mathcal{R}(T) = L^p(\Omega,\Sigma,P)$, see \cite{klw-indag, klw-exp}.

Let $E$ be a Dedekind complete Riesz space with weak order unit, say $e$. 
If $f$ is in the positive cone, $E^+:=\{ f\in E\ |\ f\ge 0\},$ of $E$ then the band generated by 
$f$ is given by
$$B_f=\{ g\in E\ |\  |g|\wedge nf\uparrow_n |g| \}.$$
Let $P_f$ be the band projection onto $B_f$, then $P_fg=P_fg^+-P_fg^-$ for $g\in E$ and
$$P_fg=\sup_{n=0,1,\dots} g\wedge nf,\quad \mbox{for}\quad g\in E^+.$$ 
In this setting, if $e$ is a weak order unit for $E$ and $P$ is a band projection onto a band $B$, then
$B$ is the principal band generated by $Pe$. 
Here we have, see \cite[Lemma 2.2]{klw-conv}, that
two non-equal elements in a Riesz space can be separated
as follows: if $m,M\in E$ with $M>m$ then there are real numbers $s<t$ so that 
$$(M-te)^+\wedge (se-m)^+>0.$$

As shown in \cite[Theorem 3.2]{klw-exp}, for $T$ a conditional expectation operator on the 
Dedekind complete Riesz space, $E$, with weak order unit $e=Te$, 
$f\in {\cal R}(T)^+$ implies $P_fT=TP_f$.  Conversely, if $Q$ is a 
 band projection on $E$ with $TQ=QT$ then $Qe\in {\cal R}(T)$ and $Q=P_{Qe}$.
For general results on bands, principal bands and band projections we refer the reader to \cite{zaanen}. 

We recall from \cite[Proposition 1.1.10]{M-N} some aspects of order convergence in $E$, 
a Dedekind complete Riesz space.
Let $(f_\alpha)$ be an order bounded net in $E$, then
$u_\alpha:=\sup\{ f_\beta : \alpha\le \beta\}$ and
$\ell_\alpha=\inf\{ f_\beta : \alpha\le \beta\}$ exist in $E$, for $\alpha$ in the index set of the net.
We denote $\lim\sup f_\alpha=\inf_\alpha u_\alpha$ and $\lim\inf f_\alpha = \sup_\alpha \ell_\alpha$.
Conversely, that both $\lim\sup f_\alpha$ and $\lim\inf f_\alpha$ exist is equivalent to requiring that $(f_\alpha)$ is order
 bounded.
Now, $(f_\alpha)$ is order convergent if and only if
$\lim\sup f_\alpha$ and $\lim\inf f_\alpha$ both exist and are equal.  In this case the common value is denoted $\lim f_\alpha$.

\begin{defs}
 Let $E$ be a Dedekind complete Riesz space with weak order unit and $T$ be a strictly positive conditional expectation on $E$.
 The space $E$ is universally complete with respect to $T$, i.e. $T$-universally complete,
 if for each increasing net $(f_\alpha)$ in $E^+$
 with $(Tf_\alpha)$ order bounded, we have that $(f_\alpha)$ is order convergent.
\end{defs}

If $E$ is a Dedekind complete Riesz space and $T$ is a strictly positive conditional 
expectation operator on $E$, then $E$ 
has a $T$-universal completion, see \cite{klw-exp}, which is 
the natural domain of $T$, denoted ${\rm dom}(T)$ in the universal completion, $E^u$, of $E$,
 also see \cite{dodds, grobler, neveu, zaanen-2}.
Here ${\rm dom}(T)=D-D$ and $Tx:=Tx^+-Tx^-$ for $x\in {\rm dom}(T)$ where
$$D=\{ x\in E^u_+ | \exists(x_\alpha)\subset E_+, x_\alpha\uparrow x, (Tx_\alpha)
   \ \mbox{order bounded in}\ E^u\},$$
and $Tx:=\sup_\alpha Tx_\alpha$, for $x\in D$ with $x_\alpha\uparrow x$,
$(x_\alpha)\subset E_+$, $(Tx_\alpha)$ order bounded in $E^u$.
It is useful to have available the following Riesz space analogues of the $L^p$ spaces as introduced in \cite{lw},
${\mathcal L}^1(T) = dom(T)$ and ${\mathcal L}^2(T)= \{x \in {\mathcal L}^1(T) \vert x^2 \in {\mathcal L}^1(T)\}.$


\newsection{$T$-conditional Independence}
The concept of $T$-conditional independence was generalized from the probability space setting to that of a Dedekind complete Riesz space, say $E$, with weak order unit, say $e$, and conditional expectation in $T$ having $Te=e$ as follows in
\cite[Definition 4.1]{klw-erg}.

\begin{defs}
 Let $E$ be a Dedekind complete Riesz space with conditional expectation $T$ and weak order unit $e=Te$.
 Let $P$ and $Q$ be band projections on $E$, we say that $P$ and $Q$ are $T$-conditionally independent with respect to $T$ if
 \begin{eqnarray}
  TPTQe=TPQe=TQTPe.\label{indep-e}
 \end{eqnarray}
 We say that two Riesz subspaces $E_1$ and $E_2$ of $E$ are $T$-conditionally independent with respect to $T$ if
 all band projections $P_i, i=1,2,$ in $E$ with  $P_ie\in E_i, i=1,2,$ are $T$-conditionally independent with respect to $T$.
\end{defs}

{\bf Example}
Consider the particular case where $E  = \mathcal{L}^1(\Omega, \mathcal{F}, P)$ is the probability space with measure $P$.  
Let $\mathcal{G}$ be a sub-$\sigma$-algebra of $\mathcal{F}$ and $T$ be the conditional expectation 
$T\cdotp = \E[\cdotp \vert \mathcal{G}]$.
The weak order units of $E$ which are invariant under $T$ are those $f \in E$ with $f>0$ almost everywhere which are 
$f \ \mathcal{G}$-measurable.  
Here the and projections are 
multiplication by characteristic functions of sets which are $\mathcal{F}$-measurable.  If we now apply Definition \ref{indep-e}, 
we have that 
$Pf = \chi_A \cdotp f$ and $Qf = \chi_B \cdotp f$, where $f$ is a weak order unit invariant under $T$.  Here $P$ and $Q$ are 
$T$-conditionally independent if 
	$$\E[ \chi_A \cdotp \E[\chi_B \cdotp f |\mathcal{G}]|\mathcal{G}] = 
		\E[ \chi_A \cdotp \chi_B \cdotp f|\mathcal{G}] = 
		\E[ \chi_B \cdotp \E[\chi_A \cdotp f |\mathcal{G}]|\mathcal{G}].$$
By the usual properties of $\E[\cdotp | \mathcal{G}]$,
$$f\E[ \chi_A |\mathcal{G}]\E[\chi_B |\mathcal{G}] = 
		f\E[ \chi_A \cdotp \chi_B  |\mathcal{G}] = 
		f\E[ \chi_B \cdotp |\mathcal{G}]\E[\chi_A |\mathcal{G}].$$
That is,
$$\E[ \chi_A |\mathcal{G}]\E[\chi_B |\mathcal{G}] = 
		\E[ \chi_A \cdotp \chi_B  |\mathcal{G}] = 
		\E[ \chi_B \cdotp |\mathcal{G}]\E[\chi_A |\mathcal{G}],$$
giving the classical definition of conditionally independent events.

Definition \ref{indep-e} is independent of the choice of the weak order unit $e$ with $e=Te$, as can be seen by the following lemma.

\begin{thm}\label{new-indep-thm}
 Let $E$ be a Dedekind complete Riesz space with conditional expectation $T$ and let $e$ 
 be a weak order unit which is invariant under $T$.
 The band projections $P$ and $Q$ in $E$ are $T$-conditionally independent 
 if, and only if,
 \begin{eqnarray}
  TPTQw=TPQw=TQTPw\quad\mbox{for all}\quad w\in \mathcal{R}(T).\label{indep-w}
 \end{eqnarray}
\end{thm}
\proof
 That (\ref{indep-w}) implies (\ref{indep-e}) is obvious.
We now show that (\ref{indep-e}) implies (\ref{indep-w}). From linearity it is sufficient
 to show that (\ref{indep-w}) holds for all $0\le w\in \mathcal{R}(T)$.\\
Consider $0\le w\in \mathcal{R}(T)$.  By Freudenthal's theorem (\cite{zaanen}), there exist 
$a_j^n \in \R$ and band projections $Q_j^n$ on $E$ such that $Q_j^n \in \mathcal{R}(T)$ and

  $$s_n=\sum_{j=0}^{n} a_j^nQ_j^ne,$$
with
 $$w=\lim_{n\to\infty} s_n.$$
 
As $e,\ Q_j^n \in \mathcal{R}(T)$,\ $Q_j^nT=TQ_j^n$.  Thus
 \begin{eqnarray}
  TPTQQ_j^ne =Q_j^nTPQe\label{indep-1}
 \end{eqnarray}
 since $Q_j^n$ commutes with all the factors in the product and therefore with the product itself.
 Again using the commutation of band projections and the fact that  $Q_j^nT=TQ_j^n$
 we obtain
 \begin{eqnarray}
  TPQQ_j^ne = Q_j^nTPQe.\label{indep-2}
 \end{eqnarray}
 Combining (\ref{indep-1}), (\ref{indep-2}) and using the linearity of $T,P$ and $Q$ gives
 \begin{eqnarray}
  TPTQ\sum_{j=0}^{n} a_j^n Q_j^ne  
       & =&TPQ\sum_{j=0}^{n} a_j^n Q_j^ne.\label{indep-3}
 \end{eqnarray}
 Since $T,P,Q$ are order continuous, taking the limit as $n\to\infty$ of (\ref{indep-3}) we obtain
 \begin{eqnarray*}
  TPTQw& =&TPQw.
 \end{eqnarray*}
 Interchanging the roles of $P$ and $Q$ gives
 \begin{eqnarray*}
  TQTPw& =&TQPw.
 \end{eqnarray*}
  As band projections commute, we have thus shown that (\ref{indep-w}) holds.
 \qed

The following corollary to the above theorem 
shows that $T$-conditional independence of the band projections $P$ and $Q$ is equivalent to
$T$-conditional independence of the closed Riesz subspaces $\left< Pe, \mathcal{R}(T)\right>$ and $\left< Qe, \mathcal{R}(T)\right>$
generated by $Pe$ and $\mathcal{R}(T)$ and by $Qe$ and $\mathcal{R}(T)$ respectively.

\begin{cor}
 Let $E$ be a Dedekind complete Riesz space with conditional expectation $T$ and let $e$ 
 be a weak order unit which is invariant under $T$.
 Let $P_i, i=1,2,$ be band projections on $E$. Then $P_i, i=1,2,$  are $T$-conditionally independent
 if and only if the closed Riesz subspaces
 $E_i=\left< P_ie, \mathcal{R}(T)\right>, i=1,2,$ are $T$-conditionally independent.
\end{cor}

\proof
 The reverse implication is obvious.  Assuming $P_i, i=1,2,$ are $T$-conditionally independent with respect to $T$
 we show that the closed Riesz subspaces $E_i, i=1,2,$ are $T$-conditionally independent with respect to $T$.
 As each element of $\mathcal{R}(T)$ is the limit of a sequence of linear combinations of band projections whose
 action on $e$ is in $\mathcal{R}(T)$ it follows that
 $E_i$ is the closure of the linear span of $$\{ P_iRe, (I-P_i)Re | R \mbox{ band projection in $E$ with } Re\in \mathcal{R}(T)\}.$$
 It thus suffices, from the linearity and continuity of band projections and conditional expectations,
 to prove that for $R_i, i=1,2,$ band projections in $E$ with
 $R_ie\in \mathcal{R}(T), i=1,2$, the band projections $P_1R_1$ and $(I-P_1)R_1$ are $T$-conditionally independent of $P_2R_2$ and $(I-P_2)R_2$.
 We will only prove that $P_1R_1$ is $T$-conditionally independent of $P_2R_2$ as the other three cases follow by similar reasoning.
 From Theorem \ref{new-indep-thm}, as $R_1e, R_2e\in \mathcal{R}(T)$,
 \begin{eqnarray*}
  TP_1TP_2R_1R_2e=TP_1P_2R_1R_2e=TP_2TP_1R_1R_2e.
 \end{eqnarray*}
 As band projections commute and since $R_iT=TR_i, i=1,2,$ we obtain 
 \begin{eqnarray*}
  TP_1R_1TP_2R_2e=TP_1R_1P_2R_2e=TP_2R_2TP_1R_1e
 \end{eqnarray*}
 giving the $T$-conditional independence of $P_iR_i, i=1,2.$ 
\qed

In the light of the above corollary, when discussing $T$-conditional independence of Riesz subspaces of $E$ with respect to $T$, we will
assume that they are closed Riesz subspaces containing $\mathcal{R}(T)$.

A Radon-Nikod\'ym-Douglas-And\^o type theorem was established in \cite{watson}. In particular, suppose
 $E$ is a $T$-universally complete Riesz space and $e=Te$ is a weak order unit, where $T$
 is a strictly positive conditional expectation operator on $E$.  
 A subset $F$ of $E$ is a closed Riesz subspace of $E$ with ${\cal R}(T)\subset F$ if and only if
 there is a unique conditional expectation $T_F$ on $E$ with
 ${\cal R}(T_F)=F$ and $TT_F=T=T_FT$. In this case $T_Ff$ for $f\in E^+$ is uniquely determined by the
 property that
 \begin{eqnarray}\label{R-N}
 TPf=TPT_Ff
 \end{eqnarray}
 for all band projections on $E$ with $Pe\in F$.
The existence and uniqueness of such conditional expectation operators forms the underlying foundation for the 
following result which characterizes independence of closed Riesz subspaces of a $T$-universally complete Riesz space 
in terms of conditional expectation operators.

\begin{thm}\label{indep}
 Let $E_1$ and $E_2$ be two closed Riesz subspaces of the $T$-universally complete Riesz space $E$ with
 strictly positive conditional expectation operator $T$ and weak order unit $e=Te$.
 Let $S$ be a conditional expectation on $E$ with $ST = T$.
 If ${\cal R}(T)\subset E_1\cap E_2$ and $T_{\left<{\cal R}(S),E_i\right>}$ is the conditional expectation having as its
range the closed Riesz space of $E$ generated by ${\cal R}(S)$ and $E_i$,  then
 the spaces $E_1$ and $E_2$ are $T$-conditionally independent with respect to $S$,  if and only if 
 $$T_iT_{\left<{\cal R}(S),E_{3-i}\right>}=T_iST_{\left<{\cal R}(S),E_{3-i}\right>} \quad i=1,2,$$
 where $T_i$ is the conditional expectation commuting with $T$ and having range $E_i$.
\end{thm}

\proof
 Let $E_1$ and $E_2$ be $T$-conditionally independent with respect to $S$, i.e.  
 for all band projections $P_i$ with $P_ie\in E_i$ for $i=1,2,$ we have
 $$SP_1SP_2e=SP_1P_2e=SP_2SP_1e.$$
 Consider the equation
 \begin{eqnarray}
 SP_1SP_2e = SP_1P_2e.\label{pf1}
 \end{eqnarray}
 Applying $T$ to both sides of the equation and using (\ref{R-N}) gives
 $$TP_1P_2e = TP_1SP_2e.$$ 
 Thus, by the Riesz space Radon-Nikod\'ym-Douglas-And\^o theorem,
 $$T_1P_2e = T_1SP_2e.$$
 Now, let $P_S$ be a band projection with $P_Se\in{\cal R}(S)$.  Applying $P_S$ and then $T$ to (\ref{pf1}) gives
 $$TP_SP_1P_2e =TP_SP_1SP_2e.$$
 As $P_Se \in {\cal R}(S)$, we have that $SP_S=P_SS$ which, together with the commutation of band projections, gives
 $$TP_1P_SP_2e =TP_1SP_SP_2e.$$
 Applying the Riesz space Radon-Nikod\'ym-Douglas-And\^o theorem now gives
  $$T_1P_SP_2e = T_1SP_SP_2e.$$
 Each element of $\left<{\cal R}(S),E_2\right>={\cal R}(T_{\left<{\cal R}(S),E_2\right>}$ can be expressed as a limit of a net of 
 linear combinations of elements of the form $P_SP_2e$ where $P_S$ and $P_2$ are respectively band projections
 with $P_Se\in {\cal R}(S)$ and $P_2e\in E_2$.  From the continuity of $T_1$
 $$T_1T_{\left<{\cal R}(S),E_2\right>} = T_1ST_{\left<{\cal R}(S),E_2\right>}.$$
 Similarly, if we consider the equation $ SP_2P_1e = SP_2SP_1e $ we have 
 $$T_2T_{\left<{\cal R}(S),E_1\right>} = T_2ST_{\left<{\cal R}(S),E_1\right>}.$$
 
Now suppose $T_iT_{\left<{\cal R}(S),E_{3-i}\right>} = T_iST_{\left<{\cal R}(S),E_{3-i}\right>}$ for all $i =1, 2$.
 Again we consider only $T_1T_{\left<{\cal R}(S),E_2\right>} = T_1ST_{\left<{\cal R}(S),E_2\right>}$. Then,
 for all $P_2e \in {\cal R}(T_2), P_Se \in {\cal R}(S)$,
 $$T_1P_SP_2e = T_1SP_SP_2e.$$
 Since $P_Se \in {\cal R}(S)$ we have
 $$T_1P_SP_2e = T_1P_SSP_2e.$$
 If we apply $P_1$, where $P_1e\in {\cal R}(T_1)$, and then $T$ to both sides of the above equality we obtain
 $$TP_1T_1P_SP_2e = TP_1T_1P_SSP_2e.$$
 Commutation of band projections, $T_1P_1=P_1T_1$ and $T=TT_1$, applied to the above equation gives
 $$TP_SP_1P_2 e = TP_SP_1SP_2e.$$
 Now from the Radon-Nikod\'ym-Douglas-And\^ o theorem in Riesz spaces we have
 $$SP_1P_2e = SP_1SP_2e.$$
 By a similar argument using $T_2T_{\left<{\cal R}(S),E_1\right>} = T_2ST_{\left<{\cal R}(S),E_1\right>}$, we have
 $$SP_2P_1e = SP_2SP_2e.$$
 Since band projections commute we get
 $$SP_1SP_2e = SP_1P_2e =SP_2SP_1e$$
 which concludes the proof.
\qed

Taking $S$ = $T$ in the above theorem, we obtain the following corollary.

\begin{cor}\label{cor-indep} 
  Let $E_1$ and $E_2$ be two closed Riesz subspaces of the $T$-universally complete Riesz space $E$ with
 strictly positive conditional expectation operator $T$ and weak order unit $e=Te$.
 If ${\cal R}(T)\subset E_1\cap E_2$, then
 the spaces $E_1$ and $E_2$ are $T$-conditionally independent,  if and only if 
 $$T_1T_2=T=T_2T_1,$$
 where $T_i$ is the conditional expectation commuting with $T$ and having range $E_i$.
\end{cor}

The following theorem is useful in the characterization of independent subspaces through conditional expectations.

\begin{cor}
 Under the same conditions as in Corollary \ref{cor-indep},
 $E_1$ and $E_2$ are $T$-conditionally independent if and only if
 \begin{eqnarray}
	T_if=Tf,\quad\mbox{for all}\quad f\in E_{3-i}, \quad i=1,2,\label{eq-indep-1}
 \end{eqnarray}
 where $T_i$ is the conditional expectation commuting with $T$ and having range $E_i$.
\end{cor}

\proof
 Observe that (\ref{eq-indep-1}) is equivalent to
	$$T_iT_{3-i}=TT_{3-i}=T.$$
 The corollary now follows directly from Corollary \ref{cor-indep}.
\qed

The above theorem can be applied to self-independence, given that the only self-independent band projections with respect to
to $T$ are those onto bands generated by elements of the range of $T$.

\begin{cor}\label{self}
  Let $E$ be a  $T$-universally complete Riesz space $E$ with
 strictly positive conditional expectation operator $T$ and weak order unit $e=Te$.
 Let $P$ be a band projection on $E$ which is self-independent with respect to $T$, then $TP=PT$ and $TPe=Pe$. 
\end{cor}

\proof
 Taking $P_1=P=P_2$ and $f=Pe$ in the above theorem, we obtain
 $TPe=T_1Pe$. But $Pe\in \mathcal{R}(T_1)$ so $TPe=Pe$, thus $Pe\in \mathcal{R}(T)$ from which it follows that $TP=PT$.
\qed

In measure theoretic probability, we can define independence of a family of $\sigma$-sub-algebras. In a similar
manner, in the Riesz space setting, we can define the independence with respect to $T$ of a family of closed Dedekind complete Riesz subspaces of $E$.\\
\\
For ease of notation, if $(E_\lambda)_{\lambda\in \Lambda}$ is a family of Riesz subspaces of $E$ we put 
$E_{\Lambda} = \left< \bigcup_{\lambda\in\Lambda_j} E_\lambda \right>$.

\begin{defs}\label{def-n}
 Let $E$ be a Dedekind complete Riesz space with conditional expectation $T$ and weak order unit $e=Te$.
 Let $E_\lambda, \lambda\in \Lambda,$ be a  family of closed Dedekind complete Riesz subspaces of $E$ having
 $\mathcal{R}(T)\subset E_\lambda$ for all $\lambda\in\Lambda$. We say that the family is $T$-conditionally independent
 if, for each pair of disjoint sets $\Lambda_1, \Lambda_2 \subset \Lambda$, we have that
 $E_{\Lambda_1}$ and $E_{\Lambda_2}$ are $T$-conditionally independent.
\end{defs}

Definition \ref{def-n} leads naturally to the definition of $T$-conditional independence for sequences in $E$, given below.

\begin{defs}\label{f-indep}
 Let $E$ be a Dedekind complete Riesz space with conditional expectation $T$ and weak order unit $e=Te$.
 We say that the sequence $(f_n)$ in $E$ is $T$-conditionally independent if
 the family $\left<\{f_n\}\cup \mathcal{R}(T)\right>,  n\in\N$ of Dedekind complete Riesz spaces is $T$-conditionally independent.
\end{defs}

\newsection{Markov Processes}

For the remainder of the paper we shall make the assumption that if $\mathcal{R}(F) \subset T$ for any closed, Dedekind complete
subspace $F$ of $E$, the conditional expectation $T_F$ onto $F$ always refers to the the unique conditional expectation 
that commutes with $T$ as is described in \ref{R-N}.\\
\\
Based on the definition of a Markov process in ${\cal L}^1$ by M. M. Rao \cite{rao} we define a Markov process in a Riesz space as follows.
\begin{defs}
 Let $T$ be a strictly positive conditional expectation on the $T$-universally complete Riesz space $E$ with weak order unit $e = Te$.
 Let $\Lambda$ be a totally ordered index set.
 A net $(X_\lambda)_{\lambda \in \Lambda}$ is a Markov process in $E$ if for any set of points $t_1 < \dots < t_n < t, \ t_i, \ t \in \Lambda$, we have  
 \begin{eqnarray}
  T_{(t_1,\dots,t_n)}Pe = T_{t_n}Pe \quad \mbox{for all} \quad Pe \in \left<{\cal R}(T), X_t\right>,  \label{MPdef}
 \end{eqnarray}
 for $P$  a band projection.  Here $T_{(t_1, t_2, \dots,t_n)}$ is the conditional expectation with range  
 $\left<{\cal R}(T), X_{t_1}, X_{t_2},\dots, X_{t_n}\right>$. 
\end{defs}

\begin{note}\label{MPnote}{\rm
 An application of Freudenthal's theorem, as in the proof of Theorem \ref{new-indep-thm}, to (\ref{MPdef}) yields that (\ref{MPdef})  is equivalent to
	$$T_{(t_1, \dots, t_n)}f = T_{t_n}f,\quad\mbox{for all} \quad f\in {\cal R}(T_t),$$
 which in turn is equivalent to
  $$T_{(t_1, \dots, t_n)}T_t = T_{t_n}T_t$$
 where $T_t$ is the conditional expectation with range 
 $\left<{\cal R}(T), X_t\right>.$ }
\end{note}

We can extend the Markov property to include the entire future, as is shown below.

\begin{lem}\label{lem-multi}
 Let $T$ be a strictly positive conditional expectation on the $T$-universally complete Riesz space $E$ with weak order unit $e = Te$.
 Let $\Lambda$ be a totally ordered index set.
 Suppose $(X_\lambda)_{\lambda \in \Lambda}$ is a Markov process in $E$. 
 If $s_m> \dots > s_1 > t > t_n > \dots > t_1$, $t_j, s_j, t \in \Lambda$ and for each $i=1,\dots,m$, $Q_i$ is a band projection with 
 $Q_ie \in \left<{\cal R}(T), X_{s_i}\right>$, then
 \begin{eqnarray}
  T_{(t_1,\dots,t_n),t}Q_1Q_2\dots Q_me = T_{t}Q_1Q_2\dots Q_me. 
  \label{MPdefFut}
 \end{eqnarray}
\end{lem}

\proof
 Under the assumptions of the lemma, if we denote $s_0=t$, from Note \ref{MPnote}
 \begin{eqnarray*}
  T_{s_j}Q_{j+1}T_{s_{j+1}}= T_{s_j}T_{s_{j+1}}Q_{j+1}= T_{(t_1,\dots,t_n,s_0,\dots,s_j)}T_{s_{j+1}}Q_{j+1}
   = T_{(t_1,\dots,t_n,s_0,\dots,s_j)}Q_{j+1}T_{s_{j+1}},
 \end{eqnarray*}
 which, if we denote $S_{s_j}=T_{(t_1,\dots,t_n,s_0,\dots,s_j)}$, gives 
 \begin{eqnarray}
  T_{s_j}Q_{j+1}T_{s_{j+1}}= S_{s_j}Q_{j+1}T_{s_{j+1}}.\label{recurrence-1}
 \end{eqnarray}
 Similarly, if we denote   $U_{s_j}=T_{(s_0,\dots,s_j)}$, then
 \begin{eqnarray}
  T_{s_j}Q_{j+1}T_{s_{j+1}}= U_{s_j}Q_{j+1}T_{s_{j+1}}.\label{recurrence-2}
 \end{eqnarray}

 Applying (\ref{recurrence-1}) recursively we obtain
  \begin{eqnarray*}
   T_{s_0}Q_1T_{s_1}Q_2T_{s_2}\dots T_{s_{m-1}}Q_me &=&  S_{s_0}Q_1T_{s_1}Q_2T_{s_2}\dots T_{s_{m-1}}Q_me\\
    &=&  S_{s_0}Q_1S_{s_1}Q_2T_{s_2}\dots T_{s_{m-1}}Q_me\\
    &=& \dots\\
    &=&  S_{s_0}Q_1S_{s_1}Q_2S_{s_2}\dots S_{s_{m-1}}Q_me.
  \end{eqnarray*}
  Here we have also used that $e=T_{s_m}e$.
  But $Q_iS_{s_j}=S_{s_j}Q_i$ and $S_{s_i}S_{s_j}=S_{s_i}$ for all $i\le j$ giving
 \begin{eqnarray*}
    S_{s_0}Q_1S_{s_1}Q_2S_{s_2}\dots S_{s_{m-1}}Q_me
    =  S_{s_0}S_{s_1}\dots S_{s_{m-1}}Q_1\dots Q_me
    =  S_{s_0}Q_1\dots Q_me.
  \end{eqnarray*}
  Combining the above two displayed equations gives
  \begin{eqnarray*}
   T_{s_0}Q_1T_{s_1}Q_2T_{s_2}\dots T_{s_{m-1}}Q_me
    =  S_{s_0}Q_1\dots Q_me.
  \end{eqnarray*}
  Similarly
  \begin{eqnarray*}
   T_{s_0}Q_1T_{s_1}Q_2T_{s_2}\dots T_{s_{m-1}}Q_me
    &=&  U_{s_0}Q_1\dots Q_me.
  \end{eqnarray*}
 Thus 
    $S_{s_0}Q_1\dots Q_me=U_{s_0}Q_1\dots Q_me$ which proves the lemma.
\qed

\begin{note}\rm
 From Freudenthal's Theorem, as in the proof of Theorem \ref{new-indep-thm}, 
 the linear span of $\left\{ Q_1\dots Q_me | Q_ie \in \left<{\cal R}(T), X_{s_i}\right>, Q_i \mbox{ band projections}, i=1,\dots,m\right\}$ is dense in 
$\left<{\cal R}(T), X_{s_1}, \dots, X_{s_m}\right>$, giving 
 \begin{eqnarray}
  T_{(t_1,\dots,t_n)}f = T_{t_n}f \label{MPdefFut2}
 \quad \mbox{for all} \quad f \in \left<{\cal R}(T), X_{s_1}, \dots, X_{s_m}\right>,  
 \end{eqnarray}
 where $s_1>s_2> \dots > s_m > t > t_n > \dots > t_1$.
\end{note}

\begin{thm}{\bf Chapman-Kolmogorov}\label{thm-KC}\\
 Let $T$ be a strictly positive conditional expectation on the $T$-universally complete Riesz space $E$ with weak order unit $e = Te$.
 Let $\Lambda$ be a totally ordered index set.
 If $(X_\lambda)_{\lambda\in\Lambda}$ is a Markov process and $u<t<n$, then
 $$T_uX=T_uT_tX,\quad \mbox{for all}\quad X\in {\cal R}(T_n),$$
 where ${\cal R}(T_u) = \left<{\cal R}(T), X_u\right>$.
\end{thm}

\proof
    We recall that $(X_\lambda)_{\lambda\in\Lambda}$ is a Markov process if for any set of points $ t_1 < \dots < t_n < t, t, t_i \in \Lambda$ one has 
	$$T_{(t_1,\dots,t_n)}X = T_{t_n}X$$ 
	where $X \in \ \left<{\cal R}(T),X_{t}\right>$. 
	Thus,
	$$ T_{(u, t)}X = T_{t}f, \quad \mbox{ for } X\in {\cal R}(T_n).$$
	Applying $T_u$ to the above equation gives
	$$T_uT_{(u, t)}X = T_uT_tX,$$
	and, thus
	$$T_uX = T_uT_tX$$
            since ${\cal R}(T_u) \subset {\cal R}(T_{(u, t)}).$
\qed

 Under the hypotheses of Theorem \ref{thm-KC}, it follows directly from the Chapman-Kolmogov Theorem and 
 Freudenthal's Theorem, as in the proof of Theorem \ref{new-indep-thm}, that if
 $(X_\lambda)_{\lambda\in\Lambda}$ is a Markov process and $u<t<n$, then
 $$T_uT_n=T_uT_tT_n.$$

It is often stated that a stochastic process
is Markov if and only if the past and future are independent given the present, see \cite[p 351]{rao}. It is clear that such independence
implies, even in the Riesz space setting, that the process is a Markov process. 
However, the non-commutation of conditional expectations onto
non-comparable closed Riesz subspaces
(or in the classical setting, the non-commutation of conditional expectations with respect to non-comparable $\sigma$-algebras),
makes the converse of the above claim more interesting.  The proof of this equivalence (part (iii) of the following theorem) relies on the fact that conditional
expectation operators are averaging operators and, in the Riesz space setting, that $E_e$ is an $f$-algebra, and is as such a commutative algebra. 
 Classical versions of the following theorem can be found in \cite{ash, getoor, rao}

\begin{thm}\label{Rao}
Let $T$ be a strictly positive conditional expectation on the $T$-universally complete Riesz space $E$ with weak order unit $e = Te$.
 Let $\Lambda$ be a totally ordered index set.
 For $(X_t)_{t\in \Lambda} \subset E$ the following are equivalent:
 \begin{enumerate}
  \item [(i)] The process, $(X_t)_{t\in \Lambda}$ is a Markov process.
  \item [(ii)] For conditional expectations
    ${\mathbb T}_u$ and $T_v$ with
    ${\cal R}({\mathbb T}_u) = \left<{\cal R}(T), X_n; n\leq u\right>$
    and ${\cal R}(T_v)= \left<\mathcal {R}(T), X_v\right>$, $u < v \text{ in } \Lambda$, we have
	$${\mathbb T}_uT_v = T_uT_v,$$
  \item [(iii)] For any $ s_m > \dots > s_1 > t > t_n > \dots > t_1$ from $\Lambda$, and
       $P,  Q$ band projections with $ Qe \in  \left< {\cal R}(T), X_{s_1}, \dots, X_{s_m}\right>$ and 
   $Pe \in  \left< \mathcal {R}(T), X_{t_1}, \dots, X_{t_n}\right>$  we have 
	$$T_tQT_tPe = T_tQPe = T_tPQe = T_tPT_tQe.$$ 
  \end{enumerate}
\end{thm}

\proof 
 {\bf (i) $\Rightarrow$ (ii)}
 Let $u<v, u,v \in \Lambda$, and $P$ be a band projection with $Pe \in \left<{\cal R}(T), X_v\right>$.
 Let $P_i$ be a band projection with $P_ie\in {\cal R}(T_{t_i}), t_1<t_2<\dots<t_n=u, n\in \N$.
 From the definition of a Markov process, for all $t_1<t_2<\dots t_n=u<t=v$ we have
 $T_{(t_1, \dots, t_n)}Pe = T_{t_n}Pe$ and $P_iT_{(t_1, \dots, t_n)}=T_{(t_1, \dots, t_n)}P_i$ thus
 $$T_{(t_1, \dots, t_n)}P_1P_2\dots P_nPe = P_1P_2\dots P_nT_{t_n}Pe.$$
 Applying $T$ to this equation gives
 \begin{eqnarray}
   TP_1P_2\dots P_nPe = TP_1P_2\dots P_nT_{t_n}Pe.\label{multi}
 \end{eqnarray}
 Note that the set of (finite) linear combinations of elements of 
 $$D=\{ P_1P_2\dots P_ne | P_i \mbox{ a band projection}, P_ie\in {\cal R}(T_{t_i}), t_1<t_2<\dots<t_n=u, n\in \N\}$$
 is dense in ${\cal R}({\mathbb T}_u)$. This together with (\ref{multi}) gives
 \begin{eqnarray}
   TQPe = TQT_{t_n}Pe\label{multi-2}
 \end{eqnarray}
 for band projections $Q$ with $Qe\in {\cal R}({\mathbb T}_u)$.
 Applying the Riesz space Radon-Nikod\'ym-Douglas-And\^ o theorem to (\ref{multi-2}) gives
 \begin{eqnarray}
   {\mathbb T}_u Pe = {\mathbb T}_uT_{t_n}Pe=T_uPe.\label{multi-3}
 \end{eqnarray}
 Now Freudenthal's theorem, as in the proof of Theorem \ref{new-indep-thm}, gives 
 \begin{eqnarray*}
   {\mathbb T}_u f= T_uf
 \end{eqnarray*}
 for $f\in {\cal R}(T_v)$, or equivalently
 \begin{eqnarray*}
   {\mathbb T}_uT_v=T_uT_v.
 \end{eqnarray*}

 {\bf (ii) $\Rightarrow$ (i) }
 Assume that for $u<v$ we have 
 \begin{eqnarray}
  {\mathbb T}_uT_v=T_uT_v.\label{mult-1}
 \end{eqnarray}
 Let $t_1 <\dots<t_n<t$.  Taking $v =t$ and $u = t_n$, we have 
 $T_{(t_1,\dots,t_n)}{\mathbb T}_u=T_{(t_1,\dots,t_n)}$ and
 $T_{(t_1,\dots,t_n)}T_u=T_u=T_{t_n}$. Thus applying $T_{(t_1,\dots,t_n)}$ to (\ref{mult-1})
 gives
 \begin{eqnarray*}
    T_{(t_1,\dots,t_n)}T_t=T_{(t_1,\dots,t_n)}{\mathbb T}_uT_v=T_{(t_1,\dots,t_n)}T_uT_v=T_{t_n}T_t.
 \end{eqnarray*}
 Applying this operator equation to $Pe$ where $P$ is a band projection with $Pe\in {\cal R}(T_t)$ gives
 that $(X_\lambda)_{\lambda\in\Lambda}$ is a Markov process.

 {\bf (i) $\Rightarrow$ (iii)}
 Let $Q$ be a band projection with
 $Qe \in  \left< {\cal R}(T), X_{s_1}, \dots, X_{s_m}\right>$  then from Lemma 
  \ref{lem-multi}
  $$T_{(t_1,\dots,t_n,t)}Qe = T_{t}Qe.$$ 
  Applying a band projection $P$ with 
   $Pe \in  \left< \mathcal {R}(T), X_{t_1}, \dots, X_{t_n}\right>$  followed by $T_t$ 
  to this equation gives
  $$T_tPQe=T_tT_{(t_1,\dots,t_n,t)}PQe =T_tPT_{(t_1,\dots,t_n,t)}Qe = T_tPT_{t}Qe.$$ 

 To prove $T_tQT_tPe = T_tQPe$, we prove $T_tQT_tPe = T_tPT_tQe$ and use the result above.
 Recall that in an $f$-algebra $Qf = Qe \cdotp f$.  Using this (the   
 commutativity of multiplication in the $f$-algebra $E_e$) and the fact that $T_t$ is an averaging operator in $E_e$ 
 we have
 \begin{eqnarray*}
  T_tQT_tPe &=& T_t((Qe)\cdot (T_tPe))\\
  &=& (T_tPe)\cdot (T_tQe) \\
  &=& (T_tQe)\cdot (T_tPe)\\
  &=& T_t((Pe\cdot (T_tQe))\\
  &=& T_tPT_tQe.
 \end{eqnarray*}
 Finally, by the commutation of band projections $T_tPQ=T_tQP$.

 {\bf (iii) $\Rightarrow$ (i)}\ 
 Suppose $ T_tPQe = T_tPT_tQe$  for all band projections $P$ and $Q$ with 
 $Qe \in \left<{\cal R}(T), X_{s_1}, \dots, X_{s_m}\right>$ and
 $Pe \in \left<{\cal R}(T), X_{t_i}, \dots, X_{t_n}\right>$. Let $R$ be a band projection 
 with $Re \in \left<{\cal R}(T), X_t\right>$,  then
 $$TRPT_{(t_1, \dots t_n, t)}Qe = TRT_{(t_1, \dots, t_n,t)}PQe= TT_{(t_1, \dots, t_n,t)}RPQe$$
 as 
 $PT_{(t_1, \dots t_n, t)} = T_{(t_1, \dots, t_n,t)}P$
 and
 $RT_{(t_1, \dots, t_n,t)}= T_{(t_1, \dots, t_n,t)}R$.
 But
 $TT_{(t_1, \dots, t_n,t)}=T=TT_t$, so
 $$TRPT_{(t_1, \dots t_n, t)}Qe =TRPQe= TT_tRPQe.$$
 Since $T_tR=RT_t$ we have
 $TT_tRPQe = TRT_tPQe$ and the hypothesis gives that
 $T_tPQe = T_tPT_tQe$ which combine to yield
 $TT_tRPQe=TRT_tPT_tQe$.
 Again appealing to the commutation of $R$ and $T_t$ and that $TT_t=T$ we have
 $$TRT_tPT_tQe= TT_tRPT_tQe=TRPT_tQe,$$
 giving
 $$TRPT_{(t_1, \dots t_n, t)}Qe =TRPT_tQe$$
 for all such $R$ and $P$.
 As the linear combinations of elements of the form $RPe$ are dense in
 $\left< {\cal R}(T), X_{t_1}, \dots, X_{t_n},X_t\right>$, we have, for all  
 $Se \in \left<{\cal R}(T), X_t, X_{t_1}, \dots, X_{t_n}\right>$, that 
	$$TST_{(t_1, \dots, t_n, t)}Qe = TST_tQe.$$
 By  (\ref{R-N}) and the unique determination of conditional expectation operators
 by their range spaces, we have that $T_{(t_1, \dots, t_n,t)}Qe =T_tQe,$
 proving the result.
\qed

\begin{note}\rm
 Proceeding in a similar manner to the proof of (i) $\Rightarrow$ (ii) in the above proof it follows that (iii) in the above theorem
 is equivalent to 
	$$\T_t\S_t =T_t= \T_t\S_t$$ 
 where $\S_t$ is the conditional expectation with range space
    ${\cal R}({\mathbb S}_u) = \left<{\cal R}(T), X_n; n\geq u\right>$.
 This shows that a process is a Markov process in a Riesz space if and only if the past and future are conditionally independent on the present.
\end{note}


\newsection{Independent Sums}

There is a natural connection between sums of independent random variables and Markov processes.  In the Riesz space case, this is illustrated by the
following theorem.

\begin{thm}\label{sum_MP}
 Let $T$ be a strictly positive conditional expectation on the $T$-universally complete Riesz space $E$ with weak order unit $e = Te$.
 Let $(f_n)$ be a sequence in $E$ which is $T$-conditionally independent then
 $$\left(\sum_{k=1}^n f_k\right)$$
 is a Markov process.
\end{thm}

\proof
 Let $S_n = \displaystyle\sum_{k=1}^{n} f_k$.
 We note that $\left<{\cal R}(T), S_1, \dots, S_n\right> \ = \  \left<{\cal R}(T), f_1, \dots, f_n\right>$.  Let $m>n$ and $P$ and $Q$ be band projections 
 with $Pe \in \  \left<{\cal R}(T), S_n\right>$ and $Qe \in \  \left<{\cal R}(T), f_{n+1},\dots,f_m\right>$.  
 Since $(f_n)$ is $T$-conditionally independent we have that $\left<{\cal R}(T), S_n\right> \subset\left<{\cal R}(T),f_1,\dots, f_n\right>$
 and $\left<{\cal R}(T), f_{n+1},\dots,f_m\right>$ are $T$-conditionally independent. Thus $P$ and $Q$ are $T$-conditionally independent with respect to $T$.

 Denote by $\T_n, T_n$ and $\S$ the conditional expectations with ranges
 $\left<{\cal R}(T),f_1,\dots, f_n\right>$, $\left<{\cal R}(T),S_n\right>$ and  $\left<{\cal R}(T),f_{n+1},\dots,f_m\right>$ respectively.
 Now from the independence of $(f_n)$ with respect to $T$ we have, by Corollary \ref{cor-indep}
 \begin{eqnarray}
  \T_n \S = T = \S\T_n.\label{sum-1}
 \end{eqnarray}
 As $Pe\in\left<{\cal R}(T), S_n\right>\  \subset \  \left<{\cal R}(T), S_1, \dots, S_n\right>$ and $\S Qe=Qe$  it follows that
 \begin{eqnarray}
  \T_nPQe = P\T_nQe= P\T_n\S Qe.\label{sum-2}
 \end{eqnarray}
 From (\ref{sum-1})
 \begin{eqnarray}
  P\T_n\S Qe=PTQe.\label{sum-3}
 \end{eqnarray}
 As ${\cal R}(T_n)\subset {\cal R}(\T_n)$, which is $T$-conditionally independent of $\S$, 
 \begin{eqnarray}
  T_n \S = T = \S T_n.\label{sum-4}
 \end{eqnarray}
 Combining (\ref{sum-3}) and (\ref{sum-4}) yields
 \begin{eqnarray}
  PTQe=PT_n\S Qe.\label{sum-5}
 \end{eqnarray}
 As noted $\S Qe=Qe$, also $T_nP=PT_n$, so 
 \begin{eqnarray}
  PT_n\S Qe=T_nPQe.\label{sum-6}
 \end{eqnarray}
 Combining (\ref{sum-2}), (\ref{sum-3}), (\ref{sum-5}) and (\ref{sum-6}) gives
 \begin{eqnarray}
  \T_nPQe = P\T_nQe= P\T_n\S Qe= PTQe= PT_n\S Qe=T_nPQe.\label{sum-7}
 \end{eqnarray}
 By Freudenthal's Theorem, as in the proof of Theorem \ref{new-indep-thm},
 the closure of the linear span of 
  $$\{ PQe | Pe \in \left<{\cal R}(T), S_n\right>, Qe \in \left<{\cal R}(T), f_{n+1},\dots,f_m\right>, P, Q \mbox{ band projections}\}$$ 
 contains ${\cal R}(T_{m})$. Thus by the order continuity of $T_n$ and $\T_n$ in (\ref{sum-7}),
 \begin{eqnarray*}
  \T_nh = T_nh
 \end{eqnarray*}
 for all $h\in \left<{\cal R}(T), S_m\right>$, proving that $(S_n)$ is a Markov process.
\qed
 
\begin{cor}\label{cormart} 
 Let $T$ be a strictly positive conditional expectation on the $T$-universally complete Riesz space $E$ with weak order unit $e = Te$.
 Let $(f_n)$ be a sequence in $E$ which is $T$-conditionally independent.
 If $Tf_i=0$ for all $i\in\N,$ then the sequence of partial sums
 $(S_n)$, where $\displaystyle{S_n=\sum_{k=1}^n  f_k},$ is a martingale with respect the filtration $(\T_n)$ where $\T_n$ is the  
 conditional expectation with range $\left<f_1,\dots,f_n,{\cal R}(T)\right>.$
\end{cor}
 
\proof
 We recall that $(F_i, \T_i)$ is a martingale if $(\T_i)_{i\in \N}$ is a filtration and 
 $F_i = \T_iF_j$, for all $i \leq j$. 
 Since ${\cal R}(\T_i) \subset {\cal R}(\T_j)$ for all $i \leq j$ we have that 
 $$\T_i\T_j = \T_i = \T_j\T_i$$
 and $(\T_n)$ is a filtration.
 Further, $f_1,\dots,f_i \in {\cal R}(\T_i)$ for all $i$ by construction of $\T_i$ giving
 $\T_iS_i=S_i$.

 If $i < j$, then from the independence of $(f_n)$ with respect to $T$ we have
 $\T_i T_j=T=T_j\T_i$ which applied to $f_j$ gives
 \begin{eqnarray}
  \T_i f_j=\T_iT_jf_j=Tf_j=0,\label{nul-sum}
 \end{eqnarray}
 Thus
 \begin{eqnarray*}
 \T_iS_j = \T_iS_i+ \sum_{k=i+1}^j \T_if_k =\T_iS_i=S_i,
 \end{eqnarray*}
 proving $(f_i,\T_i)$ a martingale.
\qed

From Corollary \ref{cormart} and \cite[Theorem 3.5]{klw-conv} we obtain the follow result regarding the
convergence of sums of independent summands. 

\begin{thm}\label{thm-conv-indep}
 Let $T$ be a strictly positive conditional expectation on the $T$-universally complete Riesz space $E$ with weak order unit $e = Te$.
 Let $(f_n)$ be a sequence in $E$ which is $T$-conditionally independent.
 If $Tf_i=0$ for all $i\in\N,$ 
  and there exists $g\in E$ such that $T\left|\sum_{i=1}^n f_i\right|\le g$ for all $n\in\N$
 then the sum
 $\sum_{k=1}^\infty f_k$ is order convergent in the sense that its sequence of partial sums is order convergent.
\end{thm}


\newsection{Brownian Motion}

The class of processes that satisfy the axioms of Brownian motion (Wiener-L\'{e}vy processes) have been generalised to
the Riesz space setting in \cite{lw}, where their martingale properties and relationship to the discrete stochastic integral were studied.
Here, in the case of $T$-universally complete spaces we show that, as in the classical $L^1$ setting, they are also Markov processes.

\begin{defs}
 Let $E$ be a Riesz space with conditional expectation $T$ and weak order unit $e = Te$.  A sequence $(f_n) \in {\mathcal L}^2(T)$
 is said to be a Brownian motion in $E$ with respect to $T$ and $e$ if 
 \begin{itemize}
  \item[(i)]  $(f_{i} - f_{i-1})$ is a $T$-conditional independent sequence where $f_0 = 0$;
  \item[(ii)] $T(f_{i} - f_{i-1}) = 0$, $i\in \N$;
  \item[(iii)] $T(f_n - f_m)^2 = \vert n - m \vert e$.
 \end{itemize}
\end{defs}

The classical definition of a Brownian motion states that the map 
$t \to f_t$ must be continuous if $\{f_t| t\in \Lambda\}$ is to be a stochastic process.
In the case where $\Lambda = \N$ this is always so.

\begin{thm}
 Let $T$ be a strictly positive conditional expectation on the $T$-universally complete Riesz space $E$ with weak order unit $e = Te$.
 Let $(f_n)$ be a Brownian motion in $E$ with respect to $T$.  Then $(f_n)$ is a Markov process.  Finally, if there exists  
 $g \in E$ such that $T|f_n|\le g$ for all $n\in \N$ (that is, the Brownian motion is $T$-bounded), then the 
 Brownian motion is order convergent.
\end{thm}
 
\proof
 Let $(f_n)$ be a Brownian motion in $E$ with respect to $T$, 
 then $(f_i - f_{i-1})_{i \in \N}$ is $T$-conditionally independent.  
 Let
 \begin{eqnarray*}
 g_1 &=& f_1 - f_0 = f_1\\
 g_2 &=& f_2 - f_1\\
 &\vdots&\\
 g_n &=& f_n - f_{n-1}.
 \end{eqnarray*}
 Here, $(g_i)_{i\in \N}$ is $T$-conditionally independent and $Tg_i=0$ for all $i\in\N$, so by
 Theorem \ref{sum_MP} the partial sums of $(g_i)$ form Markov process, i.e. $(f_n)$ is a Markov process with repsect to $T$.\\
 
  The final remark of the theorem is a direct application of Theorem \ref{thm-conv-indep}.
\qed


\end{document}